\documentclass[11pt]{amsart}

\usepackage[T1]{fontenc}
\usepackage{lmodern}

\usepackage{amsmath,amssymb,amsthm,mathtools}
\usepackage{mathrsfs}

\usepackage{microtype}

\usepackage[hidelinks]{hyperref}

\newtheorem{theorem}{Theorem}[section]
\newtheorem{proposition}[theorem]{Proposition}
\newtheorem{lemma}[theorem]{Lemma}

\newtheorem{conjecture}[theorem]{Conjecture}

\newcommand{\theoremname}{}
\newtheorem*{namedtheorem}{\theoremname}
\newenvironment{named}[1]
  {\renewcommand{\theoremname}{#1}\begin{namedtheorem}}
  {\end{namedtheorem}}

\theoremstyle{definition}
\newtheorem{definition}[theorem]{Definition}

\theoremstyle{remark}
\newtheorem{remark}[theorem]{Remark}

\title[]{Link Spectra of Lorenz-like Templates and the Horseshoe}

\author{Thiago de Paiva}
\address{Beijing International Center for Mathematical Research, Peking University, Beijing 100871, China}
\email[]{thhiagodepaiva@gmail.com}
\author{Alexander He}
\address{Sydney Mathematical Research Institute, The University of Sydney, Sydney, Australia}
\email[]{alexander.he@sydney.edu.au}

\begin{document}

\begin{abstract}
We give a complete classification of the inclusion relations between
the link spectrum of the horseshoe template $L(0,1)$ and
those of the nonnegative Lorenz-like templates $L(0,a)$. For a
template $T$, its link spectrum $\mathscr L(T)$ consists of the
link types carried by $T$, allowing parallel strands along its strips.
For integers $a\geq0$, the spectra of $L(0,a)$ and $L(0,1)$ are
comparable under inclusion if and only if $a\notin\{0,2\}$.
More precisely,
\[
\mathscr L(L(0,a))=\mathscr L(L(0,1))
\quad\Longleftrightarrow\quad a=1,
\]
while
\[
L(0,a)\prec L(0,1) \qquad (a\geq3).
\]
For $a=0$ and $a=2$, the two spectra are incomparable.

The main new containment is
\[
L(0,2n)\preceq L(0,1)\qquad(n\geq2),
\]
which goes beyond the classical parity-preserving nesting relations.
A key ingredient is a braid-theoretic criterion: every $T$-link
\[
T((r_1,s_1),\ldots,(r_k,s_k))
\]
with $s_k\geq2r_k$ is a horseshoe link. This gives a large explicit
family belonging to both the Lorenz and horseshoe spectra.

In the reverse direction, we construct a horseshoe link which is
not Lorenz, thereby disproving the Birman--Kofman conjecture that
every horseshoe link is a Lorenz link. For each integer $m\geq6$, we obtain infinitely many horseshoe links of braid index $m$ excluded from every other
nonnegative $L(0,a)$. These examples give infinitely many links
with any prescribed number of components at least three. In the
other direction, for each positive number of components, we give
infinitely many links carried by both $L(0,0)$ and $L(0,2)$ and
excluded from the horseshoe spectrum.
\end{abstract}

\maketitle

\section{Introduction}
\label{sec:introduction}

Periodic orbits of three-dimensional flows naturally form knots and
links, and their knotting and linking provide topological information
about the underlying dynamics. A fundamental tool for studying such
periodic orbit types is the theory of templates, developed by Birman
and Williams in their study of knotted periodic orbits
\cite{BirmanWilliams1983I,BirmanWilliams1983II}; see also
\cite{GhristHolmesSullivan1997}. A template is a compact branched
surface equipped with a semiflow which models the periodic orbits of
a hyperbolic invariant set. The topology of the periodic orbits can
therefore be studied through the knots and links carried by an
embedded branched surface.

For an embedded template $T\subset S^3$, let
\[
\mathscr L(T)
\]
denote the set of ambient isotopy classes of links which can be
embedded in $T$ by following its directed strips and their connections
at the branch lines. We allow parallel strands with the framing of
the strips. We refer to $\mathscr L(T)$ as the \emph{link spectrum}
of $T$. This leads to a natural comparison problem.

\medskip

\noindent
\emph{Given two templates, when is every link carried by one
also carried by the other?}

\medskip

To formulate this question, we write
\[
T_1\preceq T_2
\]
if
\[
\mathscr L(T_1)\subseteq\mathscr L(T_2),
\]
and
\[
T_1\prec T_2
\]
if
\[
\mathscr L(T_1)\subsetneq\mathscr L(T_2).
\]
Thus $\preceq$ compares the ambient isotopy types of links carried
by the templates. An embedded subtemplate inclusion implies such
a containment. To prove a containment in general, it suffices to
find, for each link carried by the first template, an ambient-isotopic
representative carried by the second.

We study this problem for the family of Lorenz-like templates
\[
L(m,n),
\]
obtained from the Lorenz template by inserting $m$ and $n$ half twists
in its two branches; see
\cite{Sullivan1993,GhristHolmesSullivan1997}. The template
\[
L(0,0)
\]
is the classical Lorenz template, while, with the conventions used
here,
\[
L(0,1)
\]
is the horseshoe template.

These templates arise from two classical models of chaotic dynamics.
The Lorenz system was introduced by Lorenz \cite{Lorenz1963}, a
geometric model for its attractor was developed by Guckenheimer and
Williams \cite{GuckenheimerWilliams1979}, and the existence of a
Lorenz attractor for the classical Lorenz equations was later
established rigorously by Tucker \cite{Tucker1999,Tucker2002}. The
horseshoe template, on the other hand, arises from a suspension of
Smale's horseshoe; Holmes and Williams constructed a template whose
periodic orbits model those of the suspended flow
\cite{HolmesWilliams1985}. These two flow models motivate the
comparison between $L(0,0)$ and $L(0,1)$.

The comparison among Lorenz-like templates goes back to Sullivan
\cite{Sullivan1993}; see also \cite{Dehornoy2011}. Sullivan's
deformation gives the nesting relation
\[
L(0,n)\preceq L(0,n-2)
\]
for every integer $n$. The corresponding subtemplate embeddings
preserve the framing of the strips
\cite[Proposition~3.2.1]{GhristHolmesSullivan1997}. These classical
containments preserve the parity of the twisting parameter.
Consequently, on the positive side
they produce two separate chains: the odd templates are naturally
related to the horseshoe template $L(0,1)$, whereas the even templates
are naturally related to the Lorenz template $L(0,0)$.

For negative twisting the picture is very different. Ghrist proved
that $L(0,n)$ is universal for every $n<0$
\cite{Ghrist1997,GhristHolmesSullivan1997}. Hence
\[
\mathscr L(L(0,n))
\]
is the set of all link types in $S^3$ whenever $n<0$. We therefore
concentrate on the nonnegative family, where the spectral containment
problem is genuinely restrictive.

Our main result gives a complete comparison between the horseshoe
template and every nonnegative Lorenz-like template.

\begin{theorem}
\label{thm:intro-comparison}
Among the nonnegative Lorenz-like templates, $L(0,a)$ is comparable
with the horseshoe template $L(0,1)$ with respect to $\preceq$ if and
only if
\[
a\notin\{0,2\}.
\]
More precisely,
\[
\mathscr L(L(0,a))
=
\mathscr L(L(0,1))
\quad\Longleftrightarrow\quad
a=1,
\]
while for every $a\geq3$,
\[
L(0,a)\prec L(0,1).
\]
For $a=0$ and $a=2$, the two templates are incomparable.
\end{theorem}

One of the main new features of
Theorem~\ref{thm:intro-comparison} is the positive even case. We prove
that
\[
L(0,2n)\preceq L(0,1)
\qquad\text{for every }n\geq2.
\]
This is somewhat unexpected, since the classical containments among
Lorenz-like templates preserve the parity of the twisting parameter.
Indeed, Sullivan's nesting relation places the positive even templates
in the chain terminating at $L(0,0)$, whereas the positive odd
templates lie in the chain terminating at $L(0,1)$. At the level
of link spectra, we obtain the additional containment
\[
L(0,2n)\preceq L(0,1)
\]
for every $n\geq2$. The proof comes from a braid-theoretic comparison
between Lorenz and horseshoe links.

The Lorenz side of this comparison is conveniently described by
$T$-links. Birman and Kofman proved that Lorenz links coincide with
$T$-links up to the convention on split trivial components
\cite{BirmanKofman2009}. Given
\[
2\leq r_1<\cdots<r_k,
\qquad
s_1,\ldots,s_k>0,
\]
the $T$-link
\[
T((r_1,s_1),\ldots,(r_k,s_k))
\]
is the closure of the positive braid
\[
(\sigma_1\cdots\sigma_{r_1-1})^{s_1}
\cdots
(\sigma_1\cdots\sigma_{r_k-1})^{s_k}.
\]
Thus $T$-links provide an explicit braid parametrisation of the
link spectrum of the Lorenz template after the split trivial components
have been set aside.

Our principal braid-theoretic containment result is the following.

\begin{theorem}
\label{thm:intro-large-T-horseshoe}
Let
\[
L=T((r_1,s_1),\ldots,(r_k,s_k))
\]
be a $T$-link. If
\[
s_k\geq2r_k,
\]
then $L$ is a horseshoe link. Consequently,
\[
L\in
\mathscr L(L(0,0))
\cap
\mathscr L(L(0,1)).
\]
\end{theorem}

An important feature of this statement is that no condition is
imposed on the preceding pairs
\[
(r_1,s_1),\ldots,(r_{k-1},s_{k-1}).
\]
Thus, after fixing arbitrary initial $T$-link data and a final width
$r_k$, every sufficiently large final twisting exponent produces a
link which is simultaneously Lorenz and horseshoe. Equivalently, the
intersection
\[
\mathscr L(L(0,0))
\cap
\mathscr L(L(0,1))
\]
contains full infinite tails in the final twisting parameter for
arbitrary fixed initial data.

This provides a substantially more flexible explicit family in the
intersection than the one exhibited in
\cite[Theorem~1.2]{dePaivaLiu2026}. In that construction, the Lorenz
and horseshoe braid descriptions are related by repeatedly using the
full-twist identity
\[
(\sigma_1\cdots\sigma_{j-1})^j
=
(\sigma_{j-1}\cdots\sigma_1)^j,
\]
and the intermediate twisting exponents are therefore required to be
multiples of the corresponding numbers of strands. By contrast,
Theorem~\ref{thm:intro-large-T-horseshoe} imposes no restriction on
the preceding twisting blocks: only the final inequality
\[
s_k\geq2r_k
\]
is required.

Theorem~\ref{thm:intro-large-T-horseshoe}, together with
Theorem~\ref{thm:even-Lorenz-T-links}, is the main ingredient in the
proof of the positive even containment
\[
L(0,2n)\preceq L(0,1)
\qquad\text{for every }n\geq2.
\]

The opposite direction of the comparison concerns a conjecture of
Birman and Kofman that every horseshoe link is a Lorenz link. The
conjecture was recorded in Kofman's lecture slides \cite{Kofman2013}; With the carrying convention of this
paper, it takes the following form.

\begin{conjecture}[Birman--Kofman]
\label{conj:Kofman}
Every horseshoe link is a Lorenz link. Equivalently,
\[
\mathscr L(L(0,1))
\subseteq
\mathscr L(L(0,0)).
\]
\end{conjecture}

We show that Conjecture~\ref{conj:Kofman} is false.

\begin{theorem}
\label{thm:intro-horseshoe-not-Lorenz}
There exists a horseshoe link which is not a Lorenz link. More
precisely, for every integer $m\geq6$, there exists a horseshoe link
$K_m$ of braid index $m$ such that
\[
K_m\notin\mathscr L(L(0,0)).
\]
In particular,
\[
\mathscr L(L(0,1))
\not\subseteq
\mathscr L(L(0,0)).
\]
\end{theorem}

Thus the inclusion in Conjecture~\ref{conj:Kofman} fails at
every braid index $m\geq6$,
and hence infinitely many pairwise distinct horseshoe links which are
not Lorenz. We show further that these same links are excluded from
every $L(0,a)$ with $a\geq0$ and $a\neq1$. By varying the twists
in their defining braids, we obtain infinitely many examples at
each fixed braid index $m\geq6$, with exactly $\lceil m/2\rceil$
components. Each lies in a chain of links with the same exclusion,
with one component added at every step.

The search for the counterexamples in
Theorem~\ref{thm:intro-horseshoe-not-Lorenz} was preceded by extensive
computational experimentation. We developed an algorithm to search for
candidate horseshoe links that might fail to admit Lorenz presentations
and carried out large-scale computations on a high-performance
computing facility. Although these computations did not directly
produce the counterexamples used here, they helped guide the subsequent
theoretical construction.

Together with the containment results described above, this leads to
a particularly simple global picture. The link spectrum of
the horseshoe template agrees with that of a nonnegative Lorenz-like
template only for $L(0,1)$ itself; it strictly contains the spectra of
all $L(0,a)$ with $a\geq3$, while $L(0,0)$ and $L(0,2)$ are
incomparable with it. Equivalently,
\[
\begin{array}{c|c}
a & \text{spectral relation with }L(0,1)\\[2mm]
\hline
0 & \text{incomparable}\\
1 & \text{equal}\\
2 & \text{incomparable}\\
a\geq3 & L(0,a)\prec L(0,1)
\end{array}
\]

In particular, the parity separation arising from the classical
nesting relations disappears in one direction at the level of
link spectra. The comparison remains asymmetric, however:
the reverse containment
\[
L(0,1)\preceq L(0,a)
\]
holds only for $a=1$.

These results therefore give a complete description of the inclusion
relations between the link spectrum of the horseshoe template
and those of all nonnegative Lorenz-like templates. Thus, within the
family $L(0,a)$, the spectral position of the horseshoe template is
completely determined.

The paper is organised as follows.
Section~\ref{sec:background} recalls the spectral containment relation
for templates, Lorenz-like templates, $T$-links, the braid description
of links carried by positive even Lorenz-like templates, and the
$W$-braid description of horseshoe links.
In Section~\ref{sec:horseshoe-T-links}, we establish the
braid-theoretic correspondence underlying
Theorem~\ref{thm:intro-large-T-horseshoe}.
Section~\ref{sec:main-template-containment} applies this result to the
positive even templates and determines precisely when
\[
L(0,a)\preceq L(0,1).
\]
In Section~\ref{sec:horseshoe-not-Lorenz}, we disprove this link
inclusion and construct, for every braid index $m\geq6$, a horseshoe
link of braid index $m$ which is not Lorenz.
Section~\ref{sec:reverse-inclusion} proves a uniform obstruction for
these links and determines precisely when
\[
L(0,1)\preceq L(0,a).
\]
Finally, Section~\ref{sec:complete-comparison} combines the two
directions to prove the comparison theorem and complete the
classification.

\section{Background and notation}
\label{sec:background}

In this section, we collect the background and notation needed for the
proofs of the main results. We first recall the notion of a template
and introduce the link-theoretic containment relation between
templates used throughout the paper. Deleting components from a
carried link gives the sublink obstruction that we use in both
directions of the comparison. We then describe $L(m,n)$ and recall
the $T^a$-link description of $L(0,a)$. For positive even $a$ this
description gives ordinary $T$-links with a bound on the final
exponent. Finally, we
recall the braid-theoretic description of horseshoe links in terms of
$W$-braids. These two braid descriptions will provide the framework
for the comparison developed in the following section.

Throughout the paper, links are considered up to ambient isotopy in
$S^3$. We use the standard generators
\[
\sigma_1,\ldots,\sigma_{m-1}
\]
of the braid group $B_m$.
We read braid factors from top to bottom and number the strands
from left to right at the top.

\subsection{Templates and link containment}
\label{subsec:template-containment}

A \emph{template} is a compact branched $2$-manifold with nonempty
boundary, equipped with a semiflow locally modelled on joining and
splitting charts. Periodic orbits of the semiflow form knots and links
in the ambient $3$-manifold.

\begin{definition}
\label{def:template-carrying}
Let $T\subset S^3$ be an embedded template. A link is \emph{carried
by} $T$ if it has an ambient-isotopic representative embedded in $T$
whose strands follow the forward direction of the strips and pass
through the branch lines using the prescribed joining and splitting
charts. Several strands may follow the same route in parallel. They
follow the framing of the strips, so a twist in a strip twists the
strands together.

The \emph{link spectrum}
\[
\mathscr L(T)
\]
is the set of ambient isotopy classes of links carried by $T$.
\end{definition}

For the Lorenz template this convention includes the parallel
push-offs used by Birman and Kofman \cite{BirmanKofman2009}.
For $L(0,a)$ it is the embedded-link convention in
\cite{dePaivaHuiRodriguez2024}.
We compare the underlying link types. The forward direction gives
each carried representative an orientation, which agrees with the
usual orientation of the braid closures used below. Genus and braid
index are taken with this orientation.

Following Sullivan's point of view of comparing templates through the
knot types carried by their periodic orbits \cite{Sullivan1993}, we
use the following link-theoretic notion of containment.

\begin{definition}
Let $T_1$ and $T_2$ be templates embedded in $S^3$. We write
\[
T_1\preceq T_2
\]
if
\[
\mathscr L(T_1)\subseteq\mathscr L(T_2).
\]
Equivalently, $T_1\preceq T_2$ if every link carried by $T_1$ can
also be carried by $T_2$.
\end{definition}

The relation $\preceq$ is reflexive and transitive. It is weaker than
containment as a subtemplate. Recall that a subtemplate $S$ of a
template $T$ is a subset $S\subset T$ which is itself a template with
the inherited strips and branch connections. A link carried by $S$
follows the same strips and connections when it is viewed in $T$.
Hence an embedded subtemplate inclusion gives
\[
T_1\subset T_2
\quad\Longrightarrow\quad
T_1\preceq T_2.
\]

Throughout this paper, containment between templates will refer to
the relation $\preceq$ unless explicitly stated otherwise.

\begin{lemma}
\label{lem:sublink-obstruction}
Let $T$ be an embedded template. Every sublink of a link carried by
$T$ is carried by $T$.
Consequently, if $K\notin\mathscr L(T)$ and $L$ contains a sublink
ambient isotopic to $K$, then $L\notin\mathscr L(T)$.
\end{lemma}

\begin{proof}
Choose a carried representative of $L$ and delete the components
outside the chosen sublink. The remaining components still follow
the same strips and the same connections at the branch lines. Thus
if $L$ is carried by $T$, so is $K$.
\end{proof}

Write $K\subset L$ when a sublink of $L$ is ambient isotopic to $K$,
allowing equality. If $K$ is a knot, it suffices for one component of
$L$ to have knot type $K$. If several template spectra exclude $K$,
the lemma applies to each of them. In particular, an excluded link
can be used to give a family of further examples by constructing
carried links that contain it. For the horseshoe, this sublink
observation also appears in \cite[Proposition~5.18]{dePaivaLiu2026}.

\subsection{Lorenz-like templates}
\label{subsec:lorenz-like-templates}

We next recall the Lorenz-like templates considered by Sullivan
\cite{Sullivan1993}. Let $m,n\in\mathbb Z$. The
\emph{$(m,n)$-Lorenz-like template}, denoted by
\[
L(m,n),
\]
is obtained from the Lorenz template by cutting along its branch line
and regluing the two strips after inserting $m$ half twists in one
strip and $n$ half twists in the other. The signs of $m$ and $n$
determine the signs of the corresponding half twists. For $L(0,a)$,
the untwisted strip means the strip corresponding to the first
parameter.

The template
\[
L(0,0)
\]
is the classical Lorenz template. We refer to the whole family
$L(m,n)$ as the family of Lorenz-like templates. Interchanging the two
strips does not change the link spectrum, so
\[
\mathscr L(L(m,n))
=
\mathscr L(L(n,m)).
\]

Links carried by $L(0,0)$ are the \emph{Lorenz links}. With the
conventions used in this paper, the template
\[
L(0,1)
\]
is the horseshoe template, and hence the horseshoe links are precisely
the elements of
\[
\mathscr L(L(0,1)).
\]

\subsection{\texorpdfstring{$T$}{T}-links and positive even Lorenz-like templates}
\label{subsec:T-links}

For $j\geq 2$, set
\[
P_j=\sigma_1\sigma_2\cdots\sigma_{j-1}.
\]
Let
\[
2\leq r_1<r_2<\cdots<r_k
\]
and let $s_1,\ldots,s_k$ be positive integers. The \emph{$T$-link}
\[
T((r_1,s_1),\ldots,(r_k,s_k))
\]
is the closure of the positive braid
\[
P_{r_1}^{s_1}
P_{r_2}^{s_2}
\cdots
P_{r_k}^{s_k}
\in B_{r_k}.
\]
Birman and Kofman proved that $T$-links coincide with Lorenz links
\cite{BirmanKofman2009}. Thus the link spectrum
\[
\mathscr L(L(0,0))
\]
consists of the unlinks and the split unions of a $T$-link with
any number of trivial components. These components can be placed
along the untwisted boundary.

We will also use the description for general $a$. Write $\Delta_r$
for the positive half twist on $r$ strands, with $\Delta_1=1$.
By \cite[Definition~3.4 and Theorem~3.5]{dePaivaHuiRodriguez2024},
the links carried by $L(0,a)$ are the closures of braids of the form
\begin{equation}
\label{eq:Ta-link-presentation}
P_{r_1}^{s_1}\cdots P_{r_k}^{s_k}\Delta_r^a,
\qquad 2\leq r_1<\cdots<r_k\leq r,
\end{equation}
together with any number of additional split trivial components.
Here $r\geq1$, all the exponents $s_i$ are positive integers, and
the product preceding $\Delta_r^a$ may be empty. The additional
trivial components are represented by unused strands outside the
$r$ strands of the displayed braid. They can be placed parallel
to the boundary of the untwisted strip. In particular, $r=1$ and
an empty prefix give an unknot to which any number of split
trivial components can be added.

We shall need a description of the links carried by positive even
Lorenz-like templates in the same language. The following formulation
is an immediate consequence of
\cite[Theorems~3.5 and~3.9]{dePaivaHuiRodriguez2024}.

\begin{theorem}
\label{thm:even-Lorenz-T-links}
Let $n$ be a positive integer. The links carried by $L(0,2n)$
are the unlinks and the split unions of any number of trivial
components with a $T$-link
\[
T((r_1,s_1),\ldots,(r_k,s_k))
\]
that has a presentation whose final pair satisfies
\[
s_k\geq nr_k.
\]
The trivial components can be placed parallel to the untwisted
boundary.
\end{theorem}

\begin{proof}
In \eqref{eq:Ta-link-presentation}, take $a=2n$. For $r\geq2$,
the final factor is $\Delta_r^{2n}=P_r^{nr}$. If the preceding
product already has a block on $r$ strands, combine it with this
factor. Otherwise append a block on $r$ strands with exponent
$nr$. In either case the final exponent is at least $nr$.
Conversely, a final block with exponent $s_k\geq nr_k$ can be
written as $P_{r_k}^{s_k-nr_k}\Delta_{r_k}^{2n}$. This gives
\eqref{eq:Ta-link-presentation}, with the block of exponent zero
omitted when $s_k=nr_k$. The case $r=1$ gives the unlinks.
\end{proof}

\subsection{Horseshoe links and \texorpdfstring{$W$}{W}-braids}
\label{subsec:horseshoe-W}

As observed above, the Lorenz-like template $L(0,1)$ is the horseshoe
template. We now recall the braid-theoretic description of its
links that will be used in the next section.

For $j\geq1$, set
\[
A_j=\sigma_j\sigma_{j-1}\cdots\sigma_1.
\]

Let $m\geq2$, and let $s_2,\ldots,s_m$ be positive integers satisfying
\[
s_m\geq m.
\]
Let
\[
m\geq w_1>w_2>\cdots>w_k\geq2
\]
be a possibly empty decreasing sequence. If this sequence is nonempty,
we also require
\[
s_m>w_1.
\]

A \emph{$W$-braid} is the positive braid on $m$ strands
\[
\begin{aligned}
\beta_W={}&
(\sigma_1\cdots\sigma_{w_1-1})
\cdots
(\sigma_1\cdots\sigma_{w_k-1})\\
&\cdot
A_1^{s_2}
A_2^{s_3}
\cdots
A_{m-1}^{s_m}.
\end{aligned}
\]
If the sequence $w_1,\ldots,w_k$ is empty, the initial product is
omitted.

The closure of $\beta_W$ is called a \emph{$W$-link} and is denoted by
\[
W((w_1,1),\ldots,(w_k,1),
(2,\overline{s_2}),\ldots,(m,\overline{s_m})).
\]
The pairs $(w_i,1)$ correspond to increasing blocks of exponent one,
whereas the barred pairs $(j,\overline{s_j})$ correspond to the
decreasing blocks
\[
A_{j-1}^{s_j}
=
(\sigma_{j-1}\sigma_{j-2}\cdots\sigma_1)^{s_j}.
\]
We also regard the trivial one-strand braid as a degenerate $W$-braid,
whose closure is the unknot.

The relevance of $W$-braids is given by the following characterization
of horseshoe links proved in \cite[Theorem~3.5]{dePaivaLiu2026}.

\begin{theorem}
\label{thm:horseshoe-W}
The horseshoe links are the unlinks and the split unions of
the closure of a $W$-braid with any number of trivial components.
These trivial components can be placed parallel to the boundary
of the untwisted strip. A $W$-braid realises the braid index of
its closure.
\end{theorem}

The split components in this statement are the unused strands in
the $T^1$-presentation \eqref{eq:Ta-link-presentation}. They may
be deleted and reinserted independently. This is the
$T^1$-presentation used in \cite[Theorem~2.4]{dePaivaLiu2026},
which is the starting point for the $W$-braid characterisation.

The two braid descriptions introduced above are therefore adapted to
the two sides of the problem considered in this paper: $T$-links
describe Lorenz links, whereas $W$-braids describe horseshoe links.
The next section compares these descriptions directly. In particular,
we shall consider $W$-braids whose initial sequence is the complete
sequence
\[
m,m-1,\ldots,2,
\]
and show that their closures admit explicit $T$-link presentations.
This comparison will yield the criterion
\[
s_k\geq2r_k
\]
for a $T$-link to be a horseshoe link and will ultimately imply the
containment
\[
L(0,2n)\preceq L(0,1),
\qquad n\geq2.
\]

\section{A large family of \texorpdfstring{$T$}{T}-links that are horseshoe links}
\label{sec:horseshoe-T-links}

In this section, we compare directly the braid descriptions of Lorenz
and horseshoe links recalled in Section~\ref{sec:background}. We first
show that a distinguished family of $W$-links admits an explicit
$T$-link representation. Reversing this correspondence then gives a
simple condition on the final pair of a $T$-link which guarantees that
the link is carried by the horseshoe template.

For $j\geq2$, set
\[
P_j=\sigma_1\sigma_2\cdots\sigma_{j-1},
\qquad
A_{j-1}=\sigma_{j-1}\sigma_{j-2}\cdots\sigma_1.
\]
We denote by $\Delta_j$ the positive half twist on $j$ strands. Thus,
\[
\Delta_j
=
A_1A_2\cdots A_{j-1}
=
P_jP_{j-1}\cdots P_2.
\]

We shall use the standard braid identity
\[
P_j\Delta_j=\Delta_jA_{j-1}.
\]
Consequently, for every integer $a\geq0$,
\[
P_j^a\Delta_j=\Delta_jA_{j-1}^a.
\]

The following elementary identity will be useful.

\begin{lemma}
\label{lem:P-to-A}
Let $a_2,\ldots,a_m$ be non-negative integers. Then
\[
P_2^{a_2}P_3^{a_3}\cdots P_m^{a_m}\Delta_m
=
A_1^{a_2+1}A_2^{a_3+1}\cdots A_{m-1}^{a_m+1}.
\]
\end{lemma}

\begin{proof}
Since
\[
\Delta_m=\Delta_{m-1}A_{m-1},
\]
we have
\[
P_m^{a_m}\Delta_m
=
\Delta_mA_{m-1}^{a_m}
=
\Delta_{m-1}A_{m-1}^{a_m+1}.
\]
Therefore,
\[
P_2^{a_2}\cdots P_m^{a_m}\Delta_m
=
P_2^{a_2}\cdots P_{m-1}^{a_{m-1}}
\Delta_{m-1}A_{m-1}^{a_m+1}.
\]
Repeating the same argument successively for
\[
P_{m-1}^{a_{m-1}}\Delta_{m-1},
\ldots,
P_2^{a_2}\Delta_2
\]
gives
\[
P_2^{a_2}\cdots P_m^{a_m}\Delta_m
=
A_1^{a_2+1}\cdots A_{m-1}^{a_m+1}.  \qedhere
\]
\end{proof}

We now establish the braid correspondence that will be used throughout
the remainder of the paper.

\begin{theorem}
\label{thm:special-W-is-T}
Let $m\geq2$. Let $L$ be a $W$-link whose initial widths are
all the integers from $m$ down to $2$, in decreasing order, and
whose decreasing blocks have exponents $s_2,\ldots,s_m$.

Let
\[
\{r_1<\cdots<r_{k-1}\}
=
\{j\in\{2,\ldots,m-1\}:s_j>1\},
\]
where this set is allowed to be empty, and set
\[
r_k=m.
\]
Define
\[
t_i=s_{r_i}-1,
\qquad 1\leq i<k,
\]
and
\[
t_k=s_m+m-1.
\]
Then
\[
L=
T((r_1,t_1),\ldots,(r_k,t_k)).
\]
Moreover,
\[
r_k=m
\qquad\text{and}\qquad
t_k\geq2r_k.
\]
\end{theorem}

\begin{proof}
By definition, $L$ is the closure of
\[
\beta=
P_mP_{m-1}\cdots P_2
A_1^{s_2}A_2^{s_3}\cdots A_{m-1}^{s_m}.
\]
Since
\[
P_mP_{m-1}\cdots P_2=\Delta_m,
\]
we have
\[
\beta=
\Delta_m
A_1^{s_2}A_2^{s_3}\cdots A_{m-1}^{s_m}.
\]

Closed braids are invariant under conjugation, and
\[
\Delta_m^{-1}\beta\Delta_m
=
A_1^{s_2}A_2^{s_3}\cdots A_{m-1}^{s_m}\Delta_m.
\]
Hence
\[
\widehat{\beta}
\cong
\widehat{
A_1^{s_2}A_2^{s_3}\cdots A_{m-1}^{s_m}\Delta_m
}.
\]

Since each $s_j$ is positive, we may apply
Lemma~\ref{lem:P-to-A} with
\[
a_j=s_j-1,
\qquad 2\leq j\leq m.
\]
It follows that
\[
A_1^{s_2}A_2^{s_3}\cdots A_{m-1}^{s_m}
=
P_2^{s_2-1}P_3^{s_3-1}\cdots
P_m^{s_m-1}\Delta_m.
\]
Therefore,
\[
\widehat{\beta}
\cong
\widehat{
P_2^{s_2-1}
P_3^{s_3-1}
\cdots
P_m^{s_m-1}
\Delta_m^2
}.
\]

The square of the positive half twist is the positive full twist, so
\[
\Delta_m^2
=
(\sigma_1\sigma_2\cdots\sigma_{m-1})^m
=
P_m^m.
\]
Hence
\[
\widehat{\beta}
\cong
\widehat{
P_2^{s_2-1}
P_3^{s_3-1}
\cdots
P_{m-1}^{s_{m-1}-1}
P_m^{s_m+m-1}
}.
\]

Whenever $s_j=1$, the factor $P_j^{s_j-1}$ is trivial and may be
omitted. By the definition of the sequence
$r_1,\ldots,r_k$, the last braid is therefore
\[
P_{r_1}^{\,s_{r_1}-1}
\cdots
P_{r_{k-1}}^{\,s_{r_{k-1}}-1}
P_{r_k}^{\,s_m+m-1}.
\]
Its closure is precisely
\[
T((r_1,t_1),\ldots,(r_k,t_k)).
\]

Finally, the initial sequence of the $W$-braid begins with
\[
w_1=m.
\]
Since it is a $W$-braid, its defining conditions give
\[
s_m>w_1=m.
\]
As $s_m$ is an integer,
\[
s_m\geq m+1.
\]
Consequently,
\[
t_k=s_m+m-1\geq2m=2r_k. \qedhere
\]
\end{proof}

The converse direction of this parameter correspondence gives the
following criterion.

\begin{named}{Theorem~\ref{thm:intro-large-T-horseshoe}}
Let
\[
L=T((r_1,s_1),\ldots,(r_k,s_k))
\]
be a $T$-link. If
\[
s_k\geq2r_k,
\]
then $L$ is a horseshoe link. Consequently,
\[
L\in
\mathscr L(L(0,0))
\cap
\mathscr L(L(0,1)).
\]
\end{named}

\begin{proof}
Set
\[
m=r_k.
\]
For each $2\leq j<m$, define
\[
a_j=
\begin{cases}
s_i, & \text{if } j=r_i \text{ for some } i<k,\\
0,   & \text{otherwise}.
\end{cases}
\]
Then $L$ is the closure of
\[
P_2^{a_2}P_3^{a_3}\cdots
P_{m-1}^{a_{m-1}}P_m^{s_k}.
\]

Define
\[
u_j=a_j+1,
\qquad 2\leq j<m,
\]
and
\[
u_m=s_k-m+1.
\]
Since
\[
s_k\geq2m,
\]
we have
\[
u_m=s_k-m+1\geq m+1>m.
\]
Therefore
\[
W=
W((m,1),(m-1,1),\ldots,(2,1),
(2,\overline{u_2}),\ldots,(m,\overline{u_m}))
\]
satisfies the defining conditions of a $W$-link.

We now apply Theorem~\ref{thm:special-W-is-T}. Since
\[
u_j>1
\quad\Longleftrightarrow\quad
a_j>0,
\]
the intermediate widths occurring in the resulting $T$-link are
precisely
\[
r_1,\ldots,r_{k-1}.
\]
Moreover, for $i<k$,
\[
u_{r_i}-1
=
a_{r_i}
=
s_i,
\]
while
\[
u_m+m-1
=
(s_k-m+1)+m-1
=
s_k.
\]
Hence the $T$-link obtained from $W$ by
Theorem~\ref{thm:special-W-is-T} is exactly
\[
T((r_1,s_1),\ldots,(r_k,s_k))
=
L.
\]

By Theorem~\ref{thm:horseshoe-W}, every $W$-link is a horseshoe link.
Thus $L$ is a horseshoe link. Since every $T$-link is a Lorenz link
\cite{BirmanKofman2009}, we conclude that
\[
L\in
\mathscr L(L(0,0))
\cap
\mathscr L(L(0,1)). \qedhere
\]
\end{proof}

For example, $T(2,5)$ is represented by the $W$-braid $P_2A_1^4$.
Here the first factor comes from the initial sequence of widths and
the second has exponent $5-2+1=4$. With two blocks, the $T$-link
$T((2,1),(3,6))$ is represented by the $W$-braid
$P_3P_2A_1^2A_2^4$. The final exponent in this example meets the
bound $s_k\geq2r_k$ with equality.

Theorem~\ref{thm:intro-large-T-horseshoe} produces a broad family of
links carried by both templates, with considerable freedom in the
parameters.

Indeed, the condition in Theorem~\ref{thm:intro-large-T-horseshoe}
depends only on the second entry of the final pair in a $T$-link
presentation. For arbitrary
\[
(r_1,s_1),\ldots,(r_{k-1},s_{k-1}),r_k,
\]
every choice
\[
s_k\geq2r_k
\]
produces a link which is simultaneously Lorenz and horseshoe.
Therefore, after fixing the initial $T$-link parameters and the final
width $r_k$, the entire tail of sufficiently large final twisting
parameters lies in
\[
\mathscr L(L(0,0))
\cap
\mathscr L(L(0,1)).
\]

It is useful to compare this family with the infinite intersection
constructed in \cite[Theorem~4.1]{dePaivaLiu2026}. In that
construction, the intermediate twisting exponents are required to be
multiples of the corresponding number of strands. The mechanism is
the identity
\[
P_j^j=A_{j-1}^j,
\]
since both sides are the positive full twist on the first $j$ strands.
Thus the common family obtained there is produced by imposing
full-twist conditions at the intermediate levels.

Theorem~\ref{thm:intro-large-T-horseshoe} imposes no such restrictions on
the preceding pairs
\[
(r_1,s_1),\ldots,(r_{k-1},s_{k-1}).
\]
The only condition is
\[
s_k\geq2r_k
\]
on the final twisting block. Hence the family obtained here has
substantially greater parametric freedom than the previously known
explicit family.

\section{Nonnegative Lorenz-like templates and the horseshoe template}
\label{sec:main-template-containment}

We first combine the description of links in positive even
Lorenz-like templates from Section~\ref{sec:background} with the
$T$-link criterion established in
Section~\ref{sec:horseshoe-T-links}. This gives the main new
containment result for the positive even family. We then combine this
result with the classical parity-preserving containments for
Lorenz-like templates and a previously known obstruction for
horseshoe links to obtain a complete characterization of the
nonnegative templates $L(0,a)$ whose link spectra are contained in
that of the horseshoe template. We finish by constructing families
of links outside the horseshoe spectrum with a prescribed number
of components.

\begin{theorem}
\label{thm:even-containment}
For every integer $n\geq2$,
\[
L(0,2n)\preceq L(0,1).
\]
Equivalently,
\[
\mathscr L(L(0,2n))
\subseteq
\mathscr L(L(0,1)).
\]
Thus every link carried by the positive even Lorenz-like template
$L(0,2n)$ is a horseshoe link.
\end{theorem}

\begin{proof}
Let
\[
L\in\mathscr L(L(0,2n)),
\qquad n\geq2.
\]

If $L$ consists only of split trivial components, then $L$ can be
realised by boundary-parallel trivial components of the horseshoe
template, and hence
\[
L\in\mathscr L(L(0,1)).
\]
We may therefore assume that $L$ has a nontrivial part.

By Theorem~\ref{thm:even-Lorenz-T-links}, $L$ is a split union of
trivial components and a $T$-link with a presentation
\[
T((r_1,s_1),\ldots,(r_k,s_k))
\]
such that
\[
s_k\geq nr_k.
\]
Since $n\geq2$, it follows that
\[
s_k\geq nr_k\geq2r_k.
\]
Therefore, by Theorem~\ref{thm:intro-large-T-horseshoe}, this $T$-link
is a horseshoe link.

The split trivial components may be reinserted along the untwisted
boundary of the horseshoe template.
Consequently,
\[
L\in\mathscr L(L(0,1)).
\]
Since $L\in\mathscr L(L(0,2n))$ was arbitrary, we conclude that
\[
\mathscr L(L(0,2n))
\subseteq
\mathscr L(L(0,1)),
\]
or equivalently,
\[
L(0,2n)\preceq L(0,1).
\]
\end{proof}

We can now determine completely which nonnegative Lorenz-like
templates of the form $L(0,a)$ have their link spectra
contained in that of the horseshoe template.

\begin{theorem}
\label{thm:main-containment}
Let $a\geq0$ be an integer. Then
\[
L(0,a)\preceq L(0,1)
\]
if and only if
\[
a\notin\{0,2\}.
\]
Equivalently,
\[
\mathscr L(L(0,a))
\subseteq
\mathscr L(L(0,1))
\]
if and only if $a\neq0,2$.
\end{theorem}

\begin{proof}
We first prove the containment when
\[
a\notin\{0,2\}.
\]

If $a=1$, the statement is immediate. Suppose that $a\geq3$ is odd.
Sullivan's parity-preserving containment, applied successively, gives
\[
L(0,a)\preceq L(0,a-2)
\preceq\cdots\preceq L(0,1)
\]
\cite{Sullivan1993,Dehornoy2011}.
Hence
\[
L(0,a)\preceq L(0,1).
\]

Now suppose that $a\geq4$ is even. Write
\[
a=2n,
\qquad n\geq2.
\]
By Theorem~\ref{thm:even-containment},
\[
L(0,a)
=
L(0,2n)
\preceq
L(0,1).
\]

It remains to show that the two excluded values $a=0$ and $a=2$ do
not satisfy the containment. In fact, the same knot provides a
counterexample in both cases. Consider the torus knot
\[
K=T(3,4)=T((3,4)).
\]

Since every $T$-link is a Lorenz link
\cite{BirmanKofman2009}, we have
\[
K\in\mathscr L(L(0,0)).
\]

Moreover, Theorem~\ref{thm:even-Lorenz-T-links}, applied with $n=1$,
shows that
\[
K\in\mathscr L(L(0,2)),
\]
since the final parameters of its $T$-link presentation satisfy
\[
4\geq3.
\]

On the other hand, $K=T(3,4)$ is the $V$-link
\[
V((3,4)).
\]
For this $V$-link, the quantity appearing in
\cite[Theorem~1.3]{dePaivaLiu2026} is
\[
E(K)
=
(4-3)(3-1)
=
2.
\]
Since $K$ is a knot and
\[
2
<
\frac{3(3-1)}{2}
=
3,
\]
\cite[Theorem~1.3]{dePaivaLiu2026} implies that $K$ is not a
$W$-link and therefore cannot be embedded in the horseshoe template.
Hence
\[
K\notin\mathscr L(L(0,1)).
\]

Thus
\[
K\in
\mathscr L(L(0,0))
\cap
\mathscr L(L(0,2)),
\qquad
K\notin\mathscr L(L(0,1)).
\]
Consequently,
\[
L(0,0)\not\preceq L(0,1)
\]
and
\[
L(0,2)\not\preceq L(0,1).
\]

Combining the two directions, we obtain
\[
L(0,a)\preceq L(0,1)
\quad\Longleftrightarrow\quad
a\notin\{0,2\}
\]
for every integer $a\geq0$.
\end{proof}

In particular, the lower bound $n\geq2$ in
Theorem~\ref{thm:even-containment} is sharp:
\[
L(0,2)\not\preceq L(0,1),
\]
whereas
\[
L(0,2n)\preceq L(0,1)
\]
for every $n\geq2$.

\begin{remark}
\label{rem:many-even-not-horseshoe}
The set
\[
\bigl(\mathscr L(L(0,0))\cap\mathscr L(L(0,2))\bigr)
\setminus\mathscr L(L(0,1))
\]
contains infinitely many distinct non-split links with each fixed
positive number of components.

Indeed, for each $p\geq3$, the torus knot $T(p,p+1)$ is carried by
$L(0,0)$ as a $T$-link and by $L(0,2)$ since $p+1\geq p$.
Its $V$-link presentation has
\[
E(T(p,p+1))=(p+1-p)(p-1)=p-1
<\frac{p(p-1)}{2}.
\]
Since $T(p,p+1)$ is a knot, the obstruction in
\cite[Theorem~1.3]{dePaivaLiu2026} shows that it is not a horseshoe
link. The braid index of $T(p,p+1)$ is $p$, since its positive
$p$-strand braid contains a full twist
\cite[Corollary~2.4]{FranksWilliams1987}. These knots are therefore
pairwise distinct.

For any fixed integer $\ell\geq1$, take $\ell$ parallel curves of
the same slope as $T(p,p+1)$ on the standard torus. Their union is
\[
T(\ell p,\ell(p+1)).
\]
This is a $T$-link whose final exponent $\ell(p+1)$ is at least
its width $\ell p$. Thus it is carried by both $L(0,0)$ and
$L(0,2)$. Every component has knot type $T(p,p+1)$, so
Lemma~\ref{lem:sublink-obstruction} excludes the whole link from the
horseshoe spectrum. For a fixed $\ell$, the component knot types
distinguish the links as $p$ varies. They are non-split, since
their defining positive braids contain a full twist on all strands
and hence every pair of components has positive linking number.

For instance, $T(6,8)$ has two components of knot type $T(3,4)$.
For a fixed $p$, increasing $\ell$ adds another parallel component
on the standard torus. The resulting link contains the preceding one
as a sublink.

One can also keep the knot $T(3,4)$ fixed while varying the link.
For $\ell\geq2$, insert $P_2^{2q}$ before the torus braid
$P_{3\ell}^{4\ell}$, with $q\geq0$. The final exponent still
satisfies the bound for $L(0,2)$. In the torus braid, the first
two strands belong to different components, since its permutation
moves positions cyclically by $4\ell$ and preserves the residue
classes modulo $\ell$. The inserted factor is pure and vanishes
on restriction to any single component. Thus every component
still has knot type $T(3,4)$. Each pair of components of the
original torus link has linking number $12$, the product of the
two slope coordinates. The inserted crossings add $q$ to the
linking number of the pair using the first two strands. The sum
of absolute pairwise linking numbers is therefore
$12\binom{\ell}{2}+q$, which distinguishes the links as $q$ varies.
\end{remark}

\section{Horseshoe links which are not Lorenz links}
\label{sec:horseshoe-not-Lorenz}

In this section, we construct a counterexample to
Conjecture~\ref{conj:Kofman} of braid index six. Two restrictions on
positive braids will be useful. The first compares crossing counts
of representatives on the same number of strands, and the second
bounds the crossings of a five-strand braid with a prescribed
linking number. We then construct two families of $W$-links
containing the six-strand example, one for each parity of the braid
index.

For $r\geq2$, set
\[
P_r=\sigma_1\sigma_2\cdots\sigma_{r-1}
\]
and
\[
A_{r-1}=\sigma_{r-1}\cdots\sigma_2\sigma_1.
\]
We write
\[
\Delta_r=A_1A_2\cdots A_{r-1}
\]
for the positive half twist on $r$ strands. Thus
\[
P_r^r=A_{r-1}^r=\Delta_r^2,
\]
where $\Delta_r^2$ is the positive full twist on $r$ strands.

We also recall the $V$-braid description of minimal braids of Lorenz
links. A $V$-braid on $p$ strands has the form
\[
D_{u_1}^{v_1}\cdots D_{u_\ell}^{v_\ell}
P_{r_1}^{s_1}\cdots P_{r_k}^{s_k}P_p^q,
\]
where
\[
D_u=
\sigma_{p-1}\sigma_{p-2}\cdots\sigma_{p-u+1},
\]
the strand parameters on each side are strictly increasing, with
$2\leq u_i\leq p$ and $2\leq r_i<p$, and
\[
q\geq p.
\]
Every non-split Lorenz link of braid index at least two has a
$V$-braid representative on its minimal number of strands
\cite[Theorem~3.9]{dePaiva2024}. The exponents in the initial
products are positive, and either product may be empty.

Throughout this section, when a link is represented as the closure of
a positive braid, its components are given the braid orientation. All
pairwise linking numbers are then nonnegative. Since changing the
orientation of a component only changes signs of the corresponding
linking numbers, the absolute values of the pairwise linking numbers
are invariants of the underlying unoriented link.

When every pair of components has nonzero linking number, an ambient
isotopy between the underlying links represented by two positive
braids identifies their braid orientations up to simultaneously
reversing all the components.
Indeed, reversing just one component of a pair changes the sign of
its linking number, which is positive in both braid orientations.
Thus any two components must agree in both orientations or be reversed
together. Simultaneous reversal preserves genus and braid index.

We begin with two elementary observations which will be used in the
proof of the basic counterexample.

\begin{lemma}
\label{lem:positive-same-strands-crossings}
Let $\beta$ and $\beta'$ be positive braid representatives of the same
oriented non-split link. If $\beta$ and $\beta'$ have the same number
of strands, then they have the same number of crossings. The same
conclusion holds if the two braids represent the same underlying
link and every pair of components has nonzero linking number.
\end{lemma}

\begin{proof}
Let $L$ be the link represented by both braids, and let $\mu$ denote
the number of components of $L$. Suppose that a positive braid
representative of $L$ has $p$ strands and $c$ crossings. By Stallings'
theorem, its canonical Seifert surface is a fibre surface
\cite{Stallings1978}. This surface is obtained from $p$ discs by
attaching $c$ positive bands, and hence has Euler characteristic
\[
\chi=p-c.
\]
Since $L$ is non-split, the fibre surface is connected. Therefore
\[
\chi=2-2g(L)-\mu,
\]
and consequently
\[
g(L)=\frac{c-p-\mu+2}{2}.
\]
The fibre surface has minimal genus, so $g(L)$ in this formula
is the genus of the oriented link. This proves the first statement.
For the second, the observation above identifies the orientations
up to simultaneous reversal, which preserves $g(L)$.
\end{proof}

The next lemma concerns braids with one component on four strands
and another on a single strand. A small linking number forces the
single strand to take part in few crossings. The order of the
increasing blocks then bounds the number of crossings of the whole
braid.

\begin{lemma}
\label{lem:five-strand-V-obstruction}
Let
\[
B=P_2^aP_3^bP_4^cP_5^d\in B_5,
\qquad
a,b,c,d\geq0,
\]
and suppose that the induced permutation of $B$ has cycle type
$(4,1)$. If the linking number of the two components of $\widehat B$
is $2$ in the braid orientation, then $B$ has at most ten crossings.
\end{lemma}

\begin{proof}
If $d=0$, the fifth strand is isolated. Since the cycle type is
$(4,1)$, this strand would be the one-strand component and its
linking number with the other component would be zero. Hence $d>0$.

Read the braid from top to bottom and let $F$ be the strand in the
fixed cycle. Every crossing of $F$ is with the other component.
The linking number is $2$ and every crossing is positive, so $F$
takes part in exactly four crossings. As $F$ returns to its starting
position, it moves to the right at two of these crossings and to
the left at the other two.

In a single factor $P_j$, the strand in position $1$ moves to
position $j$ and takes part in $j-1$ crossings. Each of the other
strands among the first $j$ moves one position to the left.
In the final factor $P_5$, the strand $F$ must move to the left,
since a move from position $1$ to position $5$ would use four
rightward crossings. Before its first rightward crossing, $F$
must reach position $1$. At most one of its leftward crossings
can occur before this point because one occurs in the final
$P_5$. It follows that $F$ starts in position $1$ or $2$.

Thus $F$ belongs to the first block. Each block leaves $F$ among
the strands in its support, and the widths of the blocks are
nondecreasing. Hence every factor $P_j$ in the braid contributes
at least one crossing of $F$.

If both rightward crossings occur in one factor, that factor is
$P_3$. The two leftward crossings occur in two further factors,
each with at most four crossings. The whole braid therefore has
at most $2+4+4=10$ crossings.

Otherwise the two rightward crossings occur in two factors $P_2$.
Between them, $F$ moves back from position $2$ to position $1$.
Since the widths are nondecreasing, this move also occurs in a
factor $P_2$. At least three of the four factors are therefore
$P_2$, and the braid has at most $3+4=7$ crossings.
\end{proof}

\subsection{The basic six-strand counterexample}

We now construct a horseshoe link whose linking numbers restrict
any Lorenz presentation to the five-strand braid considered above.
Two of its components also form a positive Hopf link, which will
be used in Section~\ref{sec:reverse-inclusion}.

Consider the $W$-link
\[
H_6=
W((2,\overline{2}),(3,\overline{1}),(4,\overline{1}),
(5,\overline{1}),(6,\overline{6})).
\]
Its defining braid is
\[
\beta_6
=
A_1^2A_2A_3A_4A_5^6.
\]
Since
\[
A_5^6=\Delta_6^2,
\]
the last factor is one positive full twist on all six strands.

By Theorem~\ref{thm:horseshoe-W}, $\beta_6$ is a minimal braid
representative of a horseshoe link. We use braid index in the
orientation of this closure. In particular,
\[
\operatorname{br}(H_6)=6.
\]

\begin{theorem}
\label{thm:H6-not-Lorenz}
The horseshoe link
\[
H_6=
W((2,\overline{2}),(3,\overline{1}),(4,\overline{1}),
(5,\overline{1}),(6,\overline{6}))
\]
is not a Lorenz link. Consequently,
\[
L(0,1)\not\preceq L(0,0).
\]
For every nonnegative integer $t$, let $H_{6,t}$ be the closure of
\[
A_1^{2t+2}A_2A_3A_4\Delta_6^2.
\]
Then $H_{6,t}$ is a horseshoe link of braid index six outside the
Lorenz spectrum.
\end{theorem}

\begin{proof}
Set
\[
\gamma=A_1^2A_2A_3A_4.
\]
The permutation induced by $\gamma$ has cycle decomposition
\[
(1\,4\,2\,5)(3)(6).
\]
Indeed, $A_1^2$ leaves the order unchanged, and $A_2$, $A_3$, and
$A_4$ move strands $3$, $4$, and $5$ to the front in turn. The
order at the bottom of the prefix is $5,4,3,1,2,6$.
Since the final full twist is pure, these are also the components of
the closure of $\beta_6$.

Let $C$ denote the component containing the four strands
\[
\{1,2,4,5\},
\]
and let $D$ and $E$ denote the one-strand components corresponding to
the strands $3$ and $6$, respectively.

The factor $A_2$ has two crossings between $C$ and $D$. Each of
$A_3$ and $A_4$ contributes one further crossing of this pair,
while $A_1^2$ contributes only crossings within $C$. Strand $6$
sits outside the prefix. Hence the prefix contributes
\[
2
\]
to $\operatorname{lk}(C,D)$ and zero to the other two pairwise
linking numbers.

The full twist $\Delta_6^2$ contains two positive crossings between
every pair of braid strands. Thus, for two distinct components
containing $a$ and $b$ braid strands, respectively, the full twist
contributes
\[
ab
\]
to their linking number. In the present case, it contributes
\[
4,\qquad 4,\qquad 1
\]
to the linking numbers of the pairs $(C,D)$, $(C,E)$, and $(D,E)$,
respectively. Therefore
\[
\operatorname{lk}(C,D)=6,
\qquad
\operatorname{lk}(C,E)=4,
\qquad
\operatorname{lk}(D,E)=1.
\]
In particular, the matrix of absolute pairwise linking numbers of
$H_6$, up to relabelling the components, is
\[
\begin{pmatrix}
0&6&4\\
6&0&1\\
4&1&0
\end{pmatrix}.
\]
Notice also that $H_6$ is non-split, since every pair of components
has nonzero linking number.

We record a further feature of this decomposition: the sublink
$D\cup E$ is the positive Hopf link. Indeed, no crossing of the prefix
$\gamma=A_1^2A_2A_3A_4$ involves both $D$ and $E$. After deleting the
four strands of $C$, the final full twist $\Delta_6^2$ restricts to
$\Delta_2^2=\sigma_1^2$ on the remaining two strands. Thus
\[
D\cup E\cong\widehat{\sigma_1^2}.
\]
This observation will also be useful in the reverse-containment
argument.

Suppose, towards a contradiction, that $H_6$ is a Lorenz link.
By the minimal $V$-braid description of Lorenz links
\cite[Theorem~3.9]{dePaiva2024}, its underlying link type admits a
minimal $V$-braid representative in the orientation of that braid.
Every pair of components has nonzero linking number. Therefore its
braid orientation agrees with that of $\beta_6$ up to simultaneous
reversal, which preserves braid index. Since
\[
\operatorname{br}(H_6)=6,
\]
this $V$-braid has six strands. Write it as
\[
\beta_V
=
D_{u_1}^{v_1}\cdots D_{u_\ell}^{v_\ell}
P_{r_1}^{s_1}\cdots P_{r_k}^{s_k}P_6^q,
\qquad
q\geq6.
\]

The braids $\beta_V$ and $\beta_6$ are positive representatives of
the same non-split link on the same number of strands. Hence, by
Lemma~\ref{lem:positive-same-strands-crossings}, they have the same
number of crossings. Since
\[
c(\beta_6)
=
2+2+3+4+30
=
41,
\]
we obtain
\[
c(\beta_V)=41.
\]

We now use the pair of components whose absolute linking number is
$1$. The factor $P_6^q$ contains the full twist
\[
P_6^6=\Delta_6^2.
\]
Suppose that these two components contain, respectively, $a$ and $b$
braid strands. The full twist alone contributes
\[
ab
\]
to their linking number with the braid orientation. Since the braid
is positive and the absolute value of their total linking number is
$1$, we have
\[
1\geq ab.
\]
Since $a,b\geq1$, it follows that
\[
a=b=1.
\]
Thus these two components are one-strand components, while the third
component contains the remaining four braid strands.

Denote the two one-strand components by $X$ and $Y$, choosing the
labels so that, with the braid orientation of $\beta_V$,
\[
\operatorname{lk}(X,C)=4,
\qquad
\operatorname{lk}(Y,C)=6,
\qquad
\operatorname{lk}(X,Y)=1,
\]
where $C$ denotes the four-strand component. This labelling is possible
because the matrix of absolute pairwise linking numbers is the one
computed above and all linking numbers of a positive braid are
nonnegative.

The full twist $\Delta_6^2$ already contributes
\[
4
\]
to $\operatorname{lk}(X,C)$ and
\[
1
\]
to $\operatorname{lk}(X,Y)$. Since the braid is positive, there can
therefore be no additional crossing involving $X$ outside this full
twist.

We claim that
\[
q=6.
\]
Indeed, if $q>6$, then, besides the full twist, the final factor
contains at least one additional copy of
\[
P_6=\sigma_1\sigma_2\sigma_3\sigma_4\sigma_5.
\]
Every braid strand participates in at least one crossing in this
factor. In particular, the strand belonging to $X$ would participate
in an additional crossing, contradicting the preceding paragraph.
Hence
\[
q=6.
\]

Consequently, the part of $\beta_V$ preceding the full twist has
\[
41-30=11
\]
crossings.

Let
\[
j\in\{1,\ldots,6\}
\]
be the position occupied by the strand belonging to $X$ before the
full twist. Since $X$ participates in no crossing before the full
twist, it remains in position $j$ throughout the prefix.

Any increasing root $P_r$ occurring in the prefix must satisfy
\[
r<j,
\]
and is therefore supported strictly to the left of the $j$-th strand.
Similarly, every decreasing root occurring in the prefix must be
supported strictly to the right of the $j$-th strand. Hence the
strands to the left and to the right of $X$ form invariant blocks
under the permutation induced by the prefix.

The cases
\[
j=3
\qquad\text{and}\qquad
j=4
\]
are impossible, since both sides then contain fewer than four strands,
whereas the component $C$ contains four braid strands.

Suppose that
\[
j=2.
\]
Then the single strand to the left of $X$ must form the component
$Y$, while the four strands to the right form $C$. Since the two
blocks have no crossings with one another before the full twist, the
full twist would give
\[
\operatorname{lk}(Y,C)=4,
\]
contrary to
\[
\operatorname{lk}(Y,C)=6.
\]
The case
\[
j=5
\]
is symmetric. Therefore
\[
j\in\{1,6\}.
\]

Suppose first that
\[
j=6.
\]
Every decreasing root involves the sixth strand, so all decreasing
factors must be absent. Hence the prefix of $\beta_V$ has the form
\[
B=P_2^aP_3^bP_4^cP_5^d
\]
for some nonnegative integers $a,b,c,d$.

Since this prefix has $11$ crossings,
\[
a+2b+3c+4d=11.
\]
Its induced permutation has cycle type $(4,1)$, corresponding to the
components $C$ and $Y$. Moreover, the full twist contributes
\[
4
\]
to $\operatorname{lk}(C,Y)$, whereas
\[
\operatorname{lk}(C,Y)=6.
\]
Thus the prefix must contribute
\[
\operatorname{lk}_{B}(C,Y)=2.
\]
This contradicts
Lemma~\ref{lem:five-strand-V-obstruction}.

Finally, suppose that
\[
j=1.
\]
Then all increasing roots are absent, and the prefix consists only of
decreasing roots supported on the last five strands. Conjugating on
those strands by the positive half twist $\Delta_5$ exchanges the
decreasing roots with the corresponding increasing roots. The widths
remain in the same order, and the closure of this five-strand braid
keeps its two components and their linking number. We therefore
obtain the same contradiction with
Lemma~\ref{lem:five-strand-V-obstruction}.

This proves that $H_6$ is outside the Lorenz spectrum.

For the extension to $H_{6,t}$, the extra factor $A_1^{2t}$ is
pure and adds crossings only within $C$. The permutation and the
pairwise linking numbers are unchanged, and $D\cup E$ is still
the same Hopf sublink. The defining braid is again a minimal
$W$-braid on six strands. If $H_{6,t}$ were Lorenz, the linking
numbers would again force the strand corresponding to $E$ to
lie at an outer position and have no crossings before the full
twist. Deleting it and conjugating by $\Delta_5$ when needed would
give a five-strand product of increasing blocks with cycle type
$(4,1)$ and linking number $2$. Its crossing count would be
$11+2t$, since the final full twist still contributes $30$
crossings. This contradicts Lemma~\ref{lem:five-strand-V-obstruction}.
\end{proof}

\subsection{Counterexamples of arbitrarily large braid index}

We now construct horseshoe links containing $H_6$ for each braid
index $m\geq6$. In each parity we identify a component whose
deletion gives the preceding example. Deleting a component of a
closed braid means removing all the strands in its permutation
cycle and relabelling the surviving strands in order. The sublink
obstruction of Lemma~\ref{lem:sublink-obstruction} then excludes the
whole family from the Lorenz spectrum.

\subsubsection*{Even braid index}

For even braid index, the family is obtained by following a crossing
on the first two strands with a half twist on the first $m-1$.
For every even integer $m\geq6$, let
\[
H_m=
W((2,\overline{2}),(3,\overline{1}),(4,\overline{1}),
\ldots,(m-1,\overline{1}),(m,\overline{m})).
\]
Its defining braid is
\[
\beta_m
=
A_1^2A_2A_3\cdots A_{m-2}A_{m-1}^{m}.
\]
Since
\[
\Delta_{m-1}=A_1A_2\cdots A_{m-2}
\]
and
\[
A_{m-1}^{m}=\Delta_m^2,
\]
we may equivalently write
\[
\beta_m
=
A_1\Delta_{m-1}\Delta_m^2.
\]

\begin{lemma}
\label{lem:even-nested-family}
For every even integer $m\geq8$, the link $H_{m-2}$ is a sublink of
$H_m$. Consequently,
\[
H_6\subset H_8\subset H_{10}\subset\cdots.
\]
\end{lemma}

\begin{proof}
Since the full twist $\Delta_m^2$ is pure, the components of $H_m$
are determined by the permutation induced by
\[
A_1\Delta_{m-1}.
\]
For even $m$, this permutation has the four-cycle
\[
(1\,\,m-2\,\,2\,\,m-1),
\]
the two-cycles
\[
(j\,\,m-j),
\qquad
3\leq j\leq\frac{m}{2}-1,
\]
and the two fixed points
\[
\frac{m}{2}
\qquad\text{and}\qquad
m.
\]
In fact, $A_1$ exchanges the first two positions and
$\Delta_{m-1}$ reverses the order of the first $m-1$. This leaves
$m/2$ fixed and pairs the other positions in this range apart
from the four-cycle. Thus there are $m/2$ components in all.
In particular,
\[
\{3,m-3\}
\]
is a two-strand component of $H_m$.

Delete this component. Geometrically, $\Delta_{m-1}$ is the positive
half twist on the first $m-1$ strands, and the deleted pair
$\{3,m-3\}$ is exchanged by its induced permutation. After deleting
these two strands and relabelling the remaining strands in their
natural order, the restriction of $\Delta_{m-1}$ is therefore the
positive half twist
\[
\Delta_{m-3}
\]
on the remaining $m-3$ strands. The factor $A_1$ is unchanged.

Likewise, deleting any two strands from the full twist
$\Delta_m^2$ and relabelling the surviving strands gives the full
twist
\[
\Delta_{m-2}^2
\]
on the remaining $m-2$ strands. Consequently, the restriction of
\[
A_1\Delta_{m-1}\Delta_m^2
\]
to the surviving strands is
\[
A_1\Delta_{m-3}\Delta_{m-2}^2.
\]
Since
\[
A_1\Delta_{m-3}\Delta_{m-2}^2
=
A_1^2A_2\cdots A_{m-4}A_{m-3}^{m-2},
\]
this is precisely the defining braid of $H_{m-2}$.

Therefore $H_{m-2}$ is a sublink of $H_m$.
\end{proof}

\subsubsection*{Odd braid index}

For odd braid index, the passage from six to seven strands adds
a one-strand component. Each subsequent step adds one component
on two strands. We first consider braid index seven. Define
\[
G_7=
W((2,\overline{1}),(3,\overline{2}),(4,\overline{1}),
(5,\overline{2}),(6,\overline{1}),(7,\overline{7})).
\]
Its defining braid is
\[
\gamma_7
=
A_1A_2^2A_3A_4^2A_5A_6^7.
\]

\begin{lemma}
\label{lem:H6-sublink-G7}
The link $H_6$ is a sublink of $G_7$.
\end{lemma}

\begin{proof}
The permutation induced by the prefix
\[
A_1A_2^2A_3A_4^2A_5
\]
has cycle decomposition
\[
(1\,5\,3\,6)(2)(4)(7).
\]
To see this, the first two blocks put the first three labels in
the order $1,3,2$. Following the next three blocks puts all seven
labels in the order $6,2,5,4,1,3,7$.
Since the final factor
\[
A_6^7=\Delta_7^2
\]
is pure, these are also the components of $G_7$.

Delete the one-strand component corresponding to the strand $2$.
The strand being deleted ends the successive blocks
$A_1,A_2^2,A_3,A_4^2,A_5$ in positions $1,3,4,1,2$. In the first
copy of $A_4$, the strand in position $5$ moves to the front and
the strand being deleted moves from position $4$ to position $5$.
The second copy moves the strand being deleted to the front.
Thus only one decreasing block remains from $A_4^2$.
Relabelling the surviving strands gives the restrictions
\[
\begin{gathered}
A_1\longmapsto 1,\qquad
A_2^2\longmapsto A_1^2,\qquad
A_3\longmapsto A_2,\\
A_4^2\longmapsto A_3,\qquad
A_5\longmapsto A_4.
\end{gathered}
\]
Moreover, the restriction of the full twist on seven strands to the
six surviving strands is the full twist on those strands, so
\[
A_6^7=\Delta_7^2
\longmapsto
\Delta_6^2=A_5^6.
\]
Hence the induced six-strand braid is
\[
A_1^2A_2A_3A_4A_5^6.
\]
This is exactly the defining braid of
\[
H_6=
W((2,\overline{2}),(3,\overline{1}),(4,\overline{1}),
(5,\overline{1}),(6,\overline{6})).
\]
Therefore
\[
H_6\subset G_7. \qedhere
\]
\end{proof}

For every odd integer $m\geq7$, define
\[
G_m=
W((2,\overline{1}),(3,\overline{2}),
(4,\overline{1}),\ldots,(m-3,\overline{1}),
(m-2,\overline{2}),(m-1,\overline{1}),
(m,\overline{m})).
\]
Its defining braid is
\[
\gamma_m
=
A_1A_2^2
\Bigl(\prod_{j=3}^{m-4} A_j\Bigr)
A_{m-3}^2A_{m-2}A_{m-1}^{m}.
\]

\begin{lemma}
\label{lem:odd-nested-family}
For every odd integer $m\geq9$, the link $G_{m-2}$ is a sublink of
$G_m$. Consequently,
\[
H_6\subset G_7\subset G_9\subset G_{11}\subset\cdots.
\]
\end{lemma}

\begin{proof}
The final factor
\[
A_{m-1}^{m}=\Delta_m^2
\]
is pure. The permutation induced by the prefix of $\gamma_m$ has the
four-cycle
\[
(1\,\,m-2\,\,3\,\,m-1),
\]
the fixed point
\[
2,
\]
the two-cycles
\[
(j\,\,m-j+1),
\qquad
4\leq j\leq\frac{m-1}{2},
\]
and the two further fixed points
\[
\frac{m+1}{2}
\qquad\text{and}\qquad
m.
\]
The first two factors put the first three labels in the order
$1,3,2$. Each factor in the product $\prod_{j=3}^{m-4}A_j$ then
moves its new label $j+1$ to the front. The factors of
$A_{m-3}^2$ next move labels $m-2$ and $2$ to the front, in that
order. Finally, $A_{m-2}$ moves label $m-1$ to the front. The order
of the labels at the bottom of the prefix is therefore
\[
m-1,\ 2,\ m-2,\ m-3,\ldots,4,\ 1,\ 3,\ m.
\]
This gives the cycles above and a total of $(m+1)/2$ components.
In particular,
\[
\{4,m-3\}
\]
is a two-strand component of $G_m$.

Delete this component and relabel the remaining $m-2$ strands in
their natural order. The first two factors involve only the first
three labels. The factor $A_3$ moves label $4$ to the front, so
it becomes trivial after deletion. In each intervening $A_j$ until
label $m-3$ reaches the front, one strand of its support has been
deleted and the label it moves survives. Thus the restrictions are
\[
A_1\longmapsto A_1,
\qquad
A_2^2\longmapsto A_2^2,
\qquad
A_3\longmapsto 1,
\]
and
\[
A_j\longmapsto A_{j-1}
\qquad
(4\leq j\leq m-5).
\]
At the end of the product, the factor $A_{m-4}$ moves the other
deleted label, $m-3$, so it also becomes trivial. In each remaining
decreasing factor the strand that moves to the front survives and
both deleted strands are in its support. This gives
\[
A_{m-4}\longmapsto 1,
\qquad
A_{m-3}^2\longmapsto A_{m-5}^2,
\qquad
A_{m-2}\longmapsto A_{m-4}.
\]
Finally, restricting the full twist to the surviving $m-2$ strands
gives
\[
A_{m-1}^{m}
=
\Delta_m^2
\longmapsto
\Delta_{m-2}^2
=
A_{m-3}^{m-2}.
\]
Consequently, the induced braid is
\[
A_1A_2^2
A_3A_4\cdots A_{m-6}
A_{m-5}^2A_{m-4}A_{m-3}^{m-2},
\]
which is precisely the defining braid of $G_{m-2}$.

Therefore
\[
G_{m-2}\subset G_m.
\]
Repeatedly applying this observation gives
\[
G_7\subset G_9\subset G_{11}\subset\cdots.
\]
Together with Lemma~\ref{lem:H6-sublink-G7}, this gives
\[
H_6\subset G_7\subset G_9\subset G_{11}\subset\cdots. \qedhere
\]
\end{proof}

We can now state the resulting infinite family.

\begin{named}{Theorem~\ref{thm:intro-horseshoe-not-Lorenz}}
For every integer $m\geq6$, there exists a horseshoe link of braid
index $m$ which is not a Lorenz link. More precisely, define
\[
K_m=
\begin{cases}
H_m, & \text{if $m$ is even},\\[2mm]
G_m, & \text{if $m$ is odd}.
\end{cases}
\]
Then $K_m$ is a horseshoe link which is not a Lorenz link and
\[
\operatorname{br}(K_m)=m.
\]
In particular, there are infinitely many pairwise distinct horseshoe
links which are not Lorenz links.
\end{named}

\begin{proof}
By construction, every $H_m$ and every $G_m$ is a $W$-link.
Therefore, by Theorem~\ref{thm:horseshoe-W}, each $K_m$ is a
horseshoe link and its defining $W$-braid is a minimal braid
representative. Hence
\[
\operatorname{br}(K_m)=m.
\]

Suppose first that $m$ is even. Repeated application of
Lemma~\ref{lem:even-nested-family} gives
\[
H_6\subset H_m.
\]

Suppose now that $m$ is odd. By
Lemmas~\ref{lem:H6-sublink-G7} and
\ref{lem:odd-nested-family}, we have
\[
H_6\subset G_7\subset G_m.
\]

Thus, in every case,
\[
H_6\subset K_m.
\]

If $K_m$ were a Lorenz link, then the sublink obstruction,
Lemma~\ref{lem:sublink-obstruction}, applied to $H_6\subset K_m$
would contradict Theorem~\ref{thm:H6-not-Lorenz}.
Therefore $K_m$ is not a Lorenz link.

Finally, each defining braid ends in a positive full twist. Thus
every pair of components has positive linking number, so isotopy
of the underlying links preserves their braid orientations up to
simultaneous reversal. Since
\[
\operatorname{br}(K_m)=m,
\]
links corresponding to distinct values of $m$ are therefore
pairwise non-isotopic.
\end{proof}

The same construction can be varied at a fixed braid index.
For each $t\geq0$, let $H_{m,t}$ be the link defined by $\beta_m$
with $A_1^2$ replaced by $A_1^{2t+2}$, and let $G_{m,t}$ be
defined by $\gamma_m$ with $A_2^2$ replaced by $A_2^{3t+2}$.
Write $K_{m,t}$ for the resulting link of either parity, so
$K_{m,0}=K_m$ and $H_{6,t}$ is the six-strand link defined above.
The added factors are powers of full twists on two and three
strands, respectively. They are pure and preserve the closure
permutations. In particular, $K_{m,t}$ has exactly
$\lceil m/2\rceil$ components. The braids are still $W$-braids
on $m$ strands, so their closures still have braid index $m$.

The two deletions in Lemmas~\ref{lem:even-nested-family}
and~\ref{lem:odd-nested-family} leave the extra factors unchanged,
since the deleted labels lie outside the first two or first three
strands, respectively. For the passage from $G_{7,t}$ to $H_{6,t}$,
strand $2$ is in the first three positions during the additional
factor $(A_2^3)^t$. Restricting that factor to the two surviving
strands gives $A_1^{2t}$, in addition to the $A_1^2$ coming from
the restriction of $A_2^2$. The other deletions in
Lemma~\ref{lem:H6-sublink-G7} are unchanged. Hence the two chains
become
\[
\begin{aligned}
H_{6,t}&\subset H_{8,t}\subset H_{10,t}\subset\cdots,\\
H_{6,t}&\subset G_{7,t}\subset G_{9,t}\subset\cdots .
\end{aligned}
\]
In each chain a step adds one component. For example,
$H_6\subset H_8\subset H_{10}$ gives links with three, four, and
five components at braid indices six, eight, and ten.

For fixed $m$, increasing $t$ by one adds two crossings to
the even braid and six to the odd braid. The number of strands
and the number of components stay the same. The genus formula in
Lemma~\ref{lem:positive-same-strands-crossings} therefore increases
the genus by one or three, respectively. The full twist in each
braid makes the link non-split and gives positive linking number
for every pair of components. Thus the genus with the braid
orientation distinguishes these underlying links as $t$ varies.

Consequently,
\[
\mathscr L(L(0,1))
\not\subseteq
\mathscr L(L(0,0)).
\]
Moreover, the failure of this containment is not confined to a
single low-complexity example: it occurs for horseshoe links of every
braid index at least six. The next section shows that the same
sublink obstruction also excludes these examples from all
$L(0,a)$ with $a\geq0$ and $a\neq1$.

\medskip

\noindent\textbf{Remark.}
The basic counterexample $H_6$ has three components. The
corresponding question for knots asks whether there is a
horseshoe knot outside the Lorenz spectrum.

\section{The reverse inclusion}
\label{sec:reverse-inclusion}

In this section, we determine when the horseshoe spectrum is
contained in the spectrum of a nonnegative Lorenz-like template.
We first use the twists in a $T^a$-braid to exclude the positive
Hopf link for $a\geq3$. Its occurrence as a sublink of $H_6$
then excludes the same six-strand example from every nonnegative
$L(0,a)$ with $a\ne1$. We use the nested families from the
preceding section to extend this result to larger links with a
prescribed number of components.

\begin{lemma}
\label{lem:Hopf-obstruction}
Let $a$ be an integer with $a\geq3$. Then the positive Hopf link
is a horseshoe link which cannot be realised in $L(0,a)$.
More generally, $L(0,a)$ carries no link with a pair of components
of absolute linking number $1$.
Consequently,
\[
L(0,1)\not\preceq L(0,a).
\]
\end{lemma}

\begin{proof}
Let $H$ denote the positive Hopf link. It is the closure of the
two-strand braid
\[
\sigma_1^2,
\]
and, in the notation of $W$-links,
\[
H=W((2,\overline{2})).
\]
Hence $H$ is a horseshoe link.

Take a link carried by $L(0,a)$ and use its braid presentation
\eqref{eq:Ta-link-presentation}. The extra trivial components have
linking number zero with every other component. For a pair of
components in the closure of the $r$-strand braid, let $u$ and $v$
be the numbers of strands belonging to the two components.
During one half twist every pair of strands crosses once. Thus each
half twist contributes exactly $uv$ positive crossings between
the components, irrespective of how the preceding half twists
have ordered the strands.

The factor $\Delta_r^a$ contributes $auv$ crossings to this pair.
The rest of the braid is positive, so the absolute linking number
satisfies
\[
2\bigl|\operatorname{lk}\bigr|\geq auv\geq3.
\]
This excludes absolute linking number $1$. The positive Hopf link
has linking number $1$ in its braid orientation. Since it is
carried by the horseshoe, we conclude that
\[
L(0,1)\not\preceq L(0,a)
\]
for every $a\geq3$.
\end{proof}

The positive Hopf sublink of $H_6$ allows us to strengthen the
counterexamples constructed in Section~\ref{sec:horseshoe-not-Lorenz}.

\begin{proposition}
\label{prop:uniform-horseshoe-obstruction}
Let $t$ be a nonnegative integer. Every horseshoe link containing
$H_{6,t}$ as a sublink is excluded from $L(0,a)$ for every
nonnegative integer $a\ne1$. In particular, for every integer $m\geq6$, the links
$K_{m,t}$ constructed in Section~\ref{sec:horseshoe-not-Lorenz}
satisfy
\[
K_{m,t}\in\mathscr L(L(0,1))\setminus
\bigcup_{\substack{a\in\mathbb Z_{\geq0}\\a\neq1}}\mathscr L(L(0,a)),
\qquad
\operatorname{br}(K_{m,t})=m.
\]
Thus the single link $H_6=K_{6,0}$ witnesses the failure of
$L(0,1)\preceq L(0,a)$ for every $a\geq0$ with $a\neq1$.
\end{proposition}

\begin{proof}
Theorem~\ref{thm:H6-not-Lorenz} shows that each $H_{6,t}$ is
carried by $L(0,1)$ and is excluded from $L(0,0)$.
Since Sullivan's nesting relation gives
\[
L(0,2)\preceq L(0,0),
\]
it is not carried by $L(0,2)$ either.

For all $t$, the components $D$ and $E$ of $H_{6,t}$ form a
positive Hopf link. For every $a\geq3$,
Lemma~\ref{lem:Hopf-obstruction} implies that this Hopf link is not
carried by $L(0,a)$. Lemma~\ref{lem:sublink-obstruction} first
excludes the six-strand link. Combining the cases gives
\[
H_{6,t}\notin\mathscr L(L(0,a))
\quad\text{for every }a\geq0\text{ with }a\neq1.
\]

Applying the lemma again gives the same exclusion for every link
containing $H_{6,t}$. In particular, the chain of sublinks for
each parity shows that $K_{m,t}$ contains $H_{6,t}$. It is a
horseshoe link of braid index $m$.
\end{proof}

We can now classify the reverse containment for every nonnegative
value of the twisting parameter.

\begin{theorem}
\label{thm:reverse-classification}
For every integer $a\geq0$,
\[
L(0,1)\preceq L(0,a)
\quad\Longleftrightarrow\quad
a=1.
\]
\end{theorem}

\begin{proof}
The case $a=1$ is immediate. For every $a\neq1$,
Proposition~\ref{prop:uniform-horseshoe-obstruction} provides a
horseshoe link $H_6$ which is not carried by $L(0,a)$. Thus
\[
L(0,1)\not\preceq L(0,a)
\]
for every $a\neq1$, which proves the assertion.
\end{proof}

\begin{remark}
\label{rem:many-horseshoe-only}
The uniform obstruction in
Proposition~\ref{prop:uniform-horseshoe-obstruction} also yields
infinitely many pairwise non-isotopic horseshoe links with any fixed
number $\ell\geq3$ of components, none of which is carried by any
$L(0,a)$ with $a\geq0$ and $a\neq1$.

For a fixed $\ell$, take the links $H_{2\ell,t}$ as $t$ varies.
They are non-split, have $\ell$ components, and are distinguished
by their genera in the braid orientation. In particular, already
at braid index six there are infinitely many three-component
examples.

More generally, for each fixed $m\geq6$ the links $K_{m,t}$ give
infinitely many examples at braid index $m$ with exactly
$\lceil m/2\rceil$ components. For any nonnegative integer $b$,
the nested families give
\[
K_{m,t}\subset K_{m+2b,t}.
\]
The larger link has exactly $b$ additional components and still
contains $H_{6,t}$. Thus the same obstruction persists as we add
any prescribed number of components along the template.

For a fixed $t$ and any $\ell\geq4$, there are also infinitely
many $\ell$-component horseshoe links that contain this same
$H_{6,t}$. In the braid for $G_{2\ell-1,t}$, replace its first
factor $A_1$ by $A_1^{2q+1}$ for $q\geq0$. The added factor is
pure, so it preserves the component count. It survives the
two-strand deletions down to seven strands, then disappears when
strand $2$ is deleted in Lemma~\ref{lem:H6-sublink-G7}. The
resulting sublink is still $H_{6,t}$. The modified braids are
$W$-braids, and the $2q$ added crossings distinguish their
closures by genus, as above.
\end{remark}

\section{Complete comparison with the horseshoe template}
\label{sec:complete-comparison}

In this section, we combine the two containment classifications.
The forward containment decides when the Lorenz-like spectrum is
contained in that of the horseshoe. The reverse containment then
distinguishes equality, strict containment, and incomparability.

\begin{named}{Theorem~\ref{thm:intro-comparison}}
Among the nonnegative Lorenz-like templates, $L(0,a)$ is comparable
with the horseshoe template $L(0,1)$ with respect to $\preceq$ if and
only if
\[
a\notin\{0,2\}.
\]
More precisely,
\[
\mathscr L(L(0,a))
=
\mathscr L(L(0,1))
\quad\Longleftrightarrow\quad
a=1,
\]
while for every $a\geq3$,
\[
L(0,a)\prec L(0,1).
\]
For $a=0$ and $a=2$, the two templates are incomparable.
\end{named}

\begin{proof}
If $a=1$, the two templates are identical and hence have the same
link spectrum.

Suppose that $a\geq3$. By Theorem~\ref{thm:main-containment},
\[
L(0,a)\preceq L(0,1).
\]
On the other hand, Theorem~\ref{thm:reverse-classification} gives
\[
L(0,1)\not\preceq L(0,a).
\]
Therefore the first containment is proper, and hence
\[
L(0,a)\prec L(0,1).
\]

Finally, suppose that $a\in\{0,2\}$. By
Theorem~\ref{thm:main-containment},
\[
L(0,a)\not\preceq L(0,1),
\]
while Theorem~\ref{thm:reverse-classification} gives
\[
L(0,1)\not\preceq L(0,a).
\]
Thus the two templates are incomparable.

The result follows.
\end{proof}

Thus, with respect to their link spectra, the nonnegative
Lorenz-like templates fall into three distinct classes relative to
the horseshoe template:
\[
a=1,
\qquad
a\geq3,
\qquad
a\in\{0,2\}.
\]
The first gives equality of spectra, the second gives proper
containment in the horseshoe spectrum, and the last consists precisely
of the two parameters for which the spectra are incomparable. The
explicit obstructions in the two directions extend to infinite
families of links with fixed numbers of components, as
shown in Remarks~\ref{rem:many-even-not-horseshoe} and
\ref{rem:many-horseshoe-only}.

\section*{Acknowledgements}

The authors would like to thank Ilya Kofman for several helpful
discussions which played an important role in the development of this
work.

The authors conceived the project and developed the original
constructions. They used a large language model to revise the English and
exposition, to check elementary crossing and linking-number
calculations, and to help formulate and check short braid arguments.
The authors take full responsibility for the mathematical statements,
proofs, references, and final text.

\bibliographystyle{amsplain}
\bibliography{Lorenz_Horseshoe_References}

\end{document}